\documentclass[11pt]{article}
\usepackage{amsmath}
\usepackage{amssymb}
\usepackage{amsthm}
\usepackage{mathrsfs}
\begin{document}
\title{Weak type estimates of intrinsic square functions on the weighted Hardy spaces}
\author{Hua Wang\,\footnote{E-mail address: wanghua@pku.edu.cn.}\\[3mm]
\footnotesize{School of Mathematical Sciences, Peking University, Beijing 100871, China}}
\date{}
\maketitle
\begin{abstract}
In this paper, by using the atomic decomposition theory of weighted Hardy spaces, we will give some weighted weak type estimates for intrinsic square functions including the Lusin area function, Littlewood-Paley $g$-function and $g^*_\lambda$-function on these spaces.\\[2mm]
\end{abstract}
\textbf{\large{1. Introduction and preliminaries}}
\par
First, let's recall some standard definitions and notations. The classical $A_p$ weight theory was first introduced by Muckenhoupt in the study of weighted
$L^p$ boundedness of Hardy-Littlewood maximal functions in [5].
Let $w$ be a nonnegative, locally integrable function defined on $\mathbb R^n$, all cubes are assumed to have their sides parallel to the coordinate axes.
We say that $w\in A_p$, $1<p<\infty$, if
$$\left(\frac1{|Q|}\int_Q w(x)\,dx\right)\left(\frac1{|Q|}\int_Q w(x)^{-\frac{1}{p-1}}\,dx\right)^{p-1}\le C \quad\mbox{for every cube}\; Q\subseteq \mathbb
R^n,$$
where $C$ is a positive constant which is independent of the choice of $Q$.\\
For the case $p=1$, $w\in A_1$, if
$$\frac1{|Q|}\int_Q w(x)\,dx\le C\,\underset{x\in Q}{\mbox{ess\,inf}}\,w(x)\quad\mbox{for every cube}\;Q\subseteq\mathbb R^n.$$
\par
It is well known that if $w\in A_p$ with $1<p<\infty$, then $w\in A_r$ for all $r>p$, and $w\in A_q$ for some $1<q<p$. We thus write $q_w\equiv\inf\{q>1:w\in A_q\}$ to denote the critical index of $w$.
\par
Given a cube $Q$ and $\lambda>0$, $\lambda Q$ denotes the cube with the same center as $Q$ whose side length is $\lambda$ times that of $Q$. $Q=Q(x_0,r)$ denotes the cube centered at $x_0$ with side length $r$. For a weight function $w$ and a measurable set $E$, we set the weighted measure $w(E)=\int_E w(x)\,dx$.
\par
We give the following result that will often be used in the sequel.
\newtheorem*{lemmaA}{Lemma A}
\begin{lemmaA}[{[2]}]
Let $w\in A_p$, $p\ge1$. Then, for any cube $Q$, there exists an absolute constant $C>0$ such that
$$w(2Q)\le C w(Q).$$
In general, for any $\lambda>1$, we have
$$w(\lambda Q)\le C\lambda^{np}w(Q),$$
where $C$ does not depend on $Q$ nor on $\lambda$.
\end{lemmaA}
\par
Given a Muckenhoupt's weight function $w$ on $\mathbb R^n$, for $0<p<\infty$, we denote by $L^p_w(\mathbb R^n)$ the space of all functions satisfying
$$\|f\|_{L^p_w(\mathbb R^n)}=\left(\int_{\mathbb R^n}|f(x)|^pw(x)\,dx\right)^{1/p}<\infty.$$
We also denote by $WL^p_w(\mathbb R^n)$ the weak weighted $L^p$ space which is formed by all functions satisfying
$$\|f\|_{WL^p_w(\mathbb R^n)}=\sup_{\lambda>0}\lambda\cdot w\big(\{x\in\mathbb R^n:|f(x)|>\lambda\}\big)^{1/p}<\infty.$$
\par
For any $0<p<\infty$, the weighted Hardy spaces $H^p_w(\mathbb R^n)$ can be defined in terms of maximal functions.
Let $\varphi$ be a function in $\mathscr S(\mathbb R^n)$ satisfying $\int_{\mathbb R^n}\varphi(x)\,dx=1$.
Set
$$\varphi_t(x)=t^{-n}\varphi(x/t),\quad t>0,\;x\in\mathbb R^n.$$
We will define the maximal function $M_\varphi f(x)$ by
$$M_\varphi f(x)=\sup_{t>0}|f*\varphi_t(x)|.$$
Then $H^p_w(\mathbb R^n)$ consists of those tempered distributions $f\in\mathscr S'(\mathbb R^n)$ for which
$M_\varphi f\in L^p_w(\mathbb R^n)$ with $\|f\|_{H^p_w}=\|M_\varphi f\|_{L^p_w}$. For every $1<p<\infty$, as in the unweighted case, we have $L^p_w(\mathbb R^n)=H^p_w(\mathbb R^n)$.
\par
The real-variable theory of weighted Hardy spaces have been studied by many authors. In 1979, Garcia-Cuerva studied the atomic decomposition and the dual spaces of $H^p_w$ for $0<p\le1$. In 2002, Lee and Lin gave the molecular characterization of $H^p_w$ for $0<p\le1$, they also obtained the $H^p_w(\mathbb R)$, $\frac12<p\le1$ boundedness of the Hilbert transform and the $H^p_w(\mathbb R^n)$, $\frac n{n+1}<p\le1$ boundedness of the Riesz transforms. For the results mentioned above, we refer the readers to [1,3,6] for further details.
\par
In this article, we will use Garcia-Cuerva's atomic decomposition theory for weighted Hardy spaces in [1,6]. We characterize weighted Hardy spaces in terms of atoms in the following way.
\par
Let $0<p\le1\le q\le\infty$ and $p\ne q$ such that $w\in A_q$ with critical index $q_w$. Set [\,$\cdot$\,] the greatest integer function. For $s\in \mathbb Z_+$ satisfying $s\ge[n({q_w}/p-1)],$ a real-valued function $a(x)$ is called ($p,q,s$)-atom centered at $x_0$ with respect to $w$(or $w$-($p,q,s$)-atom centered at $x_0$) if the following conditions are satisfied:
\par (a) $a\in L^q_w(\mathbb R^n)$ and is supported in a cube $Q$ centered at $x_0$,
\par (b) $\|a\|_{L^q_w}\le w(Q)^{1/q-1/p}$,
\par (c) $\int_{\mathbb R^n}a(x)x^\alpha\,dx=0$ for every multi-index $\alpha$ with $|\alpha|\le s$.
\newtheorem*{theoremb}{Theorem B}
\begin{theoremb}
Let $0<p\le1\le q\le\infty$ and $p\ne q$ such that $w\in A_q$ with critical index $q_w$. For each $f\in H^p_w(\mathbb R^n)$, there exist a sequence \{$a_j$\} of $w$-$(p,q,[n(q_w/p-1)])$-atoms and a sequence \{$\lambda_j$\} of real numbers with $\sum_j|\lambda_j|^p\le C\|f\|^p_{H^p_w}$ such that $f=\sum_j\lambda_j a_j$ both in the sense of distributions and in the $H^p_w$ norm.
\end{theoremb}
\noindent\textbf{\large{2. The intrinsic square functions and our main results}}
\par
The intrinsic square functions were first introduced by Wilson in [8] and [9],
the so-called intrinsic square functions are defined as follows.
For $0<\alpha\le1$, let ${\mathcal C}_\alpha$ be the family of functions $\varphi$ defined on $\mathbb R^n$ such that $\varphi$ has support containing in $\{x\in\mathbb R^n: |x|\le1\}$, $\int_{\mathbb R^n}\varphi(x)\,dx=0$ and for all $x, x'\in \mathbb R^n$,
$$|\varphi(x)-\varphi(x')|\le|x-x'|^\alpha.$$
For $(y,t)\in {\mathbb R}^{n+1}_+ =\mathbb R^n\times(0,\infty)$ and $f\in L^1_{{loc}}(\mathbb R^n)$, we set
$$A_\alpha(f)(y,t)=\sup_{\varphi\in{\mathcal C}_\alpha}|f*\varphi_t(y)|.$$
Then we define the intrinsic square function of $f$(of order $\alpha$) by the formula
$$S_\alpha(f)(x)=\left(\iint_{\Gamma(x)}\Big(A_\alpha(f)(y,t)\Big)^2\frac{dydt}{t^{n+1}}\right)^{1/2},$$
where $\Gamma(x)$ denotes the usual cone of aperture one:
$$\Gamma(x)=\{(y,t)\in{\mathbb R}^{n+1}_+:|x-y|<t\}.$$
We can also define varying-aperture versions of $S_\alpha(f)$ by the formula
$$S_{\alpha,\beta}(f)(x)=\left(\iint_{\Gamma_\beta(x)}\Big(A_\alpha(f)(y,t)\Big)^2\frac{dydt}{t^{n+1}}\right)^{1/2},$$
where $\Gamma_\beta(x)$ is the usual cone of aperture $\beta>0$:
$$\Gamma_\beta(x)=\{(y,t)\in{\mathbb R}^{n+1}_+:|x-y|<\beta t\}.$$
The intrinsic Littlewood-Paley $g$-function(could be viewed as ``zero-aperture" version of $S_\alpha(f)$) and the intrinsic $g^*_\lambda$-function(could be viewed as ``infinite aperture" version of $S_\alpha(f)$) will be defined respectively by
$$g_\alpha(f)(x)=\left(\int_0^\infty\Big(A_\alpha(f)(x,t)\Big)^2\frac{dt}{t}\right)^{1/2}$$
and
$$g^*_{\lambda,\alpha}(f)(x)=\left(\iint_{{\mathbb R}^{n+1}_+}\left(\frac t{t+|x-y|}\right)^{\lambda n}\Big(A_\alpha(f)(y,t)\Big)^2\frac{dydt}{t^{n+1}}\right)^{1/2}.$$
\par
In [9], Wilson proved that the intrinsic square functions are bounded operators on the weighted Lebesgue spaces $L^p_w(\mathbb R^n)$ for $1<p<\infty$, namely, he showed the following result.
\newtheorem*{Theoremc}{Theorem C}
\begin{Theoremc}
Let $w\in A_p$, $1<p<\infty$ and $0<\alpha\le1$. Then there exists a positive constant $C>0$ such that$$\|S_\alpha(f)\|_{L^p_w}\le C \|f\|_{L^p_w}.$$
\end{Theoremc}
In [7], the authors considered some boundedness properties of intrinsic square functions on the weighted Hardy spaces $H^p_w(\mathbb R^n)$ for $0<p<1$. Moreover, they gave the intrinsic square function characterizations of weighted Hardy spaces $H^p_w(\mathbb R^n)$ for $0<p<1$. As a continuation of [7], the main purpose of this paper is to study their weak type estimates on these spaces.
\par
In order to state our theorems, we need to introduce the Lipschitz space $Lip(\alpha,1,0)$ for $0<\alpha\le1$.
Set $b_Q=\frac{1}{|Q|}\int_Q b(x)\,dx$.
$$Lip(\alpha,1,0)=\{b\in L_{loc}(\mathbb R^n):\|b\|_{Lip(\alpha,1,0)}<\infty\},$$
where
$$\|b\|_{Lip(\alpha,1,0)}=\sup_Q\frac{1}{|Q|^{1+\alpha/n}}\int_Q|b(y)-b_Q|\,dy$$
and the supremum is taken over all cubes $Q$ in $\mathbb R^n$.
\par
Our main results are stated as follows.
\newtheorem{theorem1}{Theorem}
\begin{theorem1}
Let $0<\alpha<1$, $p=n/{(n+\alpha)}$ and $w\in A_1$. Suppose that $f\in \big(Lip(\alpha,1,0)\big)^*$, then there exists a constant $C>0$ independent of $f$ such that
$$\|g_\alpha(f)\|_{WL^p_w}\le C\|f\|_{H^p_w}.$$
\end{theorem1}
\begin{theorem1}
Let $0<\alpha<1$, $p=n/{(n+\alpha)}$ and $w\in A_1$. Suppose that $f\in \big(Lip(\alpha,1,0)\big)^*$, then there exists a constant $C>0$ independent of $f$ such that
$$\|S_\alpha(f)\|_{WL^p_w}\le C\|f\|_{H^p_w}.$$
\end{theorem1}
\begin{theorem1}
Let $0<\alpha<1$, $p=n/{(n+\alpha)}$, $w\in A_1$ and $\lambda>3+(2\alpha)/n$. Suppose that $f\in \big(Lip(\alpha,1,0)\big)^*$, then there exists a constant $C$ independent of $f$ such that
$$\|g^*_{\lambda,\alpha}(f)\|_{WL^p_w}\le C\|f\|_{H^p_w}.$$
\end{theorem1}
\noindent\textbf{Remark.} Clearly, if for every $t>0$, $\varphi_t\in \mathcal C _\alpha$, then we have $\varphi_t\in Lip(\alpha,1,0)$. Thus the intrinsic square functions are well defined for tempered distributions in $\big(Lip(\alpha,1,0)\big)^*$.
\par
Throughout this article, we will use $C$ to denote a positive constant, which is independent of the main parameters and not necessarily the same at each occurrence.\\
\textbf{\large{3. Proofs of Theorems 1 and 2}}
\par
By adopting the same method given in [4, page 123], we can prove the following superposition principle on the weighted weak type estimates.
\newtheorem*{lemma3.1}{Lemma 3.1}
\begin{lemma3.1}
Let $w\in A_1$ and $0<p<1$.
If a sequence of measurable functions $\{f_j\}$ satisfy
$$\|f_j\|_{WL^p_w}\le1 \quad\mbox{for all}\;\, j\in\mathbb Z$$
and
$$\sum_{j\in\mathbb Z}|\lambda_j|^p\le1,$$
then we have
$$\Big\|\sum_{j\in\mathbb Z}\lambda_j f_j\Big\|^p_{WL^p_w}\le\frac{2-p}{1-p}.$$
\end{lemma3.1}
\begin{proof}[Proof of Theorem 1]
First we observe that for $w\in A_1$ and $p=n/{(n+\alpha)}$, then $[n({q_w}/p-1)]=[\alpha]=0$. By Theorem B and Lemma 3.1, it suffices to show that for any $w$-$(p,q,0)$-atom $a$, there exists a constant $C>0$ independent of $a$ such that $\|g_\alpha(a)\|_{WL^p_w}\le C$.
\par
Let $a$ be a $w$-$(p,q,0)$-atom with $supp\, a\subseteq Q=Q(x_0,r)$, and let $Q^*=2\sqrt nQ$. For any given $\lambda>0$, we write
\begin{equation*}
\begin{split}
&\lambda^p\cdot w(\{x\in\mathbb R^n:|g_\alpha(a)(x)|>\lambda\})\\
\le&\lambda^p\cdot w(\{x\in Q^*:|g_\alpha(a)(x)|>\lambda\})+\lambda^p\cdot w(\{x\in(Q^*)^c:|g_\alpha(a)(x)|>\lambda\})\\
=&I_1+I_2.
\end{split}
\end{equation*}
Since $w\in A_1$, then $w\in A_q$ for $1<q<\infty$. Applying Chebyshev's inequality, H\"older's inequality, Lemma A and Theorem C, we thus have
\begin{equation}
\begin{split}
I_1&\le\int_{Q^*}|g_\alpha(a)(x)|^pw(x)\,dx\\
&\le\left(\int_{Q^*}|g_\alpha(a)(x)|^qw(x)\,dx\right)^{p/q}\left(\int_{Q^*}w(x)\,dx\right)^{1-p/q}\\
&\le \|g_\alpha(a)\|^p_{L^q_w}w(Q^*)^{1-p/q}\\
&\le C\cdot\|a\|^p_{L^q_w}w(Q)^{1-p/q}\\
&\le C.
\end{split}
\end{equation}
\par
We now turn to estimate $I_2$.
For any $\varphi\in{\mathcal C}_\alpha$, $0<\alpha<1$, by the vanishing moment condition of atom $a$, we have
\begin{equation}
\begin{split}
\big|a*\varphi_t(x)\big|&=\left|\int_Q\big(\varphi_t(x-y)-\varphi_t(x-x_0)\big)a(y)\,dy\right|\\
&\le\int_Q\frac{|y-x_0|^\alpha}{t^{n+\alpha}}|a(y)|\,dy\\
&\le C\cdot\frac{r^\alpha}{t^{n+\alpha}}\int_Q|a(y)|\,dy.
\end{split}
\end{equation}
For any fixed $q>1$, we denote the conjugate exponent of $q$ by $q'=q/{(q-1)}$. H\"{o}lder's inequality and the $A_q$ condition yield
\begin{equation}
\begin{split}
\int_Q|a(y)|\,dy&\le\left(\int_Q|a(y)|^qw(y)\,dy\right)^{1/q}\left(\int_Q w(y)^{-1/(q-1)}\,dy\right)^{1/q'}\\
&\le C\cdot\|a\|_{L^q_w}\left(\frac{|Q|^q}{w(Q)}\right)^{1/q}\\
&\le C\cdot\frac{|Q|}{w(Q)^{1/p}}.
\end{split}
\end{equation}
Observe that $supp \,\varphi\subseteq\{x\in\mathbb R^n:|x|\le1\}$, then for any $y\in Q$, $x\in(Q^*)^c$, we have $t\ge|x-y|\ge|x-x_0|-|y-x_0|\ge\frac{|x-x_0|}{2}$. Substituting the above inequality (3) into (2), we thus obtain
\begin{equation*}
\begin{split}
\big|g_\alpha(a)(x)\big|^2&=\int_0^\infty\Big(\sup_{\varphi\in{\mathcal C}_\alpha}\big|a*\varphi_t(x)\big|\Big)^2\frac{dt}{t}\\
&\le C\left(\frac{r^{n+\alpha}}{w(Q)^{1/p}}\right)^2\int_{\frac{|x-x_0|}{2}}^\infty\frac{dt}{t^{2(n+\alpha)+1}}\\
&\le C\left(\frac{r^{n+\alpha}}{w(Q)^{1/p}}\right)^2\frac{1}{|x-x_0|^{2n+2\alpha}}\\
&\le C\left(\frac{1}{w(Q)^{1/p}}\right)^2.
\end{split}
\end{equation*}
Set $Q^*_0=Q$, $Q^*_1=Q^*$ and $Q^*_k=(Q^*_{k-1})^*, k=2,3,\ldots.$ Following the same lines as above, we can also show that for any $x\in(Q^*_k)^c$, then
\begin{equation*}
\big|g_\alpha(a)(x)\big|\le C\cdot\frac{1}{w(Q^*_{k-1})^{1/p}}\quad k=1,2,\ldots.
\end{equation*}
If $\{x\in(Q^*)^c:|g_\alpha(a)(x)|>\lambda\}=\O$, then the inequality
$$I_2\le C$$
holds trivially.
\\
If $\{x\in(Q^*)^c:|g_\alpha(a)(x)|>\lambda\}\neq\O$. For $p=n/{(n+\alpha)}$, it is easy to verify that $$\lim_{k\to\infty}\frac{1}{w(Q^*_k)^{1/p}}=0.$$
Then for any fixed $\lambda>0$, we are able to find a maximal positive integer $K$ such that $$\lambda<C\cdot\frac{1}{w(Q^*_K)^{1/p}}.$$
Therefore
\begin{equation}
\begin{split}
I_2&\le \lambda^p\cdot\sum_{k=1}^K w\big(\{x\in Q^*_{k+1}\backslash Q^*_k:|g_\alpha(a)(x)|>\lambda\}\big)\\
&\le C\cdot\frac{1}{w(Q^*_K)}\sum_{k=1}^K w(Q^*_{k+1})\\
&\le C.
\end{split}
\end{equation}
Combining the above inequality (4) with (1) and taking the supremum over all $\lambda>0$, we complete the proof of Theorem 1.
\end{proof}
\begin{proof}[Proof of Theorem 2]
The proof is almost the same. We only point out the main differences. For any given $\lambda>0$, we write
\begin{equation*}
\begin{split}
&\lambda^p\cdot w(\{x\in\mathbb R^n:|S_\alpha(a)(x)|>\lambda\})\\
\le&\lambda^p\cdot w(\{x\in Q^*:|S_\alpha(a)(x)|>\lambda\})+\lambda^p\cdot w(\{x\in(Q^*)^c:|S_\alpha(a)(x)|>\lambda\})\\
=&J_1+J_2.
\end{split}
\end{equation*}
Using the same arguments as in the proof of Theorem 1, we can prove
$$J_1\le C.$$
\par
To estimate $J_2$, we note that $z\in Q$, $x\in (Q^*)^c$, then $|z-x_0|\le\frac{|x-x_0|}2$. Furthermore, when $|x-y|<t$ and $|y-z|<t$, then we deduce
$$2t>|x-z|\ge|x-x_0|-|z-x_0|\ge\frac{|x-x_0|}2.$$
By using the inequalities (2) and (3), we thus obtain
\begin{equation*}
\begin{split}
|S_\alpha(a)(x)|^2&=\iint_{\Gamma(x)}\Big(\sup_{\varphi\in{\mathcal C}_\alpha}|a*\varphi_t(y)|\Big)^2\frac{dydt}{t^{n+1}}\\
&\le C\left(\frac{r^{n+\alpha}}{w(Q)^{1/p}}\right)^2\int_{\frac{|x-x_0|}{4}}^\infty\int_{|y-x|<t}\frac{dydt}{t^{2n+2\alpha+n+1}}\\
&\le C\left(\frac{r^{n+\alpha}}{w(Q)^{1/p}}\right)^2\frac{1}{|x-x_0|^{2n+2\alpha}}\\
&\le C\left(\frac{1}{w(Q)^{1/p}}\right)^2,
\end{split}
\end{equation*}
which is equivalent to
$$|S_\alpha(a)(x)|\le C\cdot\frac{1}{w(Q)^{1/p}}.$$
The rest of the proof is exactly the same as that of Theorem 1, we can get
$$J_2\le C.$$
This completes the proof of Theorem 2.
\end{proof}
\noindent\textbf{\large{4. Proof of Theorem 3}}
\par
Before proving our main theorem, we need to establish the following lemma.
\newtheorem*{lemma4.1}{Lemma 4.1}
\begin{lemma4.1}
Let $w\in A_1$ and $0<\alpha\le1$. Then for every $\lambda>1$, we have
$$\|g^*_{\lambda,\alpha}(a)\|_{L^2_w}\le C\|a\|_{L^2_w}.$$
\end{lemma4.1}
\begin{proof}
From the definition, we readily see that
\begin{align}
\big(g^*_{\lambda,\alpha}(a)(x)\big)^2&=\iint_{\mathbb R^{n+1}_+}\left(\frac{t}{t+|x-y|}\right)^{\lambda n}(A_\alpha(a)(y,t))^2\frac{dydt}{t^{n+1}}\notag\\
&=\int_0^\infty\int_{|x-y|<t}\left(\frac{t}{t+|x-y|}\right)^{\lambda n}(A_\alpha(a)(y,t))^2\frac{dydt}{t^{n+1}}\notag\\
&+\sum_{k=1}^\infty\int_0^\infty\int_{2^{k-1}t\le|x-y|<2^kt}\left(\frac{t}{t+|x-y|}\right)^{\lambda n}(A_\alpha(a)(y,t))^2\frac{dydt}{t^{n+1}}\notag\\
&\le C\bigg[S_\alpha(a)(x)^2+\sum_{k=1}^\infty 2^{-k\lambda n}S_{\alpha,2^k}(a)(x)^2\bigg].
\end{align}
We are now going to estimate $\int_{\mathbb R^n}|S_{\alpha,2^k}(a)(x)|^2w(x)\,dx$ for $k=1,2,\ldots.$
Fubini theorem and Lemma A imply 
\begin{align}
&\int_{\mathbb R^n}|S_{\alpha,2^k}(a)(x)|^2w(x)\,dx\notag\\
=&\int_{\mathbb R^n}\bigg(\int_{{\mathbb R}^{n+1}_+}\Big(A_\alpha(a)(y,t)\Big)^2\chi_{|x-y|<2^k t}\frac{dydt}{t^{n+1}}\bigg)w(x)\,dx\notag\\
=&\int_{{\mathbb R}^{n+1}_+}\Big(\int_{|x-y|<2^k t}w(x)\,dx\Big)\Big(A_\alpha(a)(y,t)\Big)^2\frac{dydt}{t^{n+1}}\notag\\
\le& \,C\cdot2^{kn}\int_{{\mathbb R}^{n+1}_+}\Big(\int_{|x-y|<t}w(x)\,dx\Big)\Big(A_\alpha(a)(y,t)\Big)^2\frac{dydt}{t^{n+1}}\\
=&\,C\cdot 2^{kn}\|S_\alpha(a)\|_{L^2_w}^2\notag.
\end{align}
Since $w\in A_1$, then $w\in A_2$. Therefore, by using Theorem C and the above inequality (6), we thus obtain
\begin{equation*}
\begin{split}
&\|g^*_{\lambda,\alpha}(a)\|^2_{L^2_w}\\
\le& \,C\bigg(\int_{\mathbb R^n}|S_\alpha(a)(x)|^2w(x)\,dx+\sum_{k=1}^\infty 2^{-k\lambda n}\int_{\mathbb R^n}|S_{\alpha,2^k}(a)(x)|^2w(x)\,dx\bigg)\\
\le& \,C\bigg(\|S_\alpha(a)\|^2_{L^2_w}+\sum_{k=1}^\infty 2^{-k\lambda n}\cdot2^{kn}\|S_\alpha(a)\|_{L^2_w}^2\bigg)\\
\end{split}
\end{equation*}
\begin{equation*}
\begin{split}
\le& \,C\cdot\|a\|^2_{L^2_w}\Big(1+\sum_{k=1}^\infty 2^{-k\lambda n}\cdot2^{kn}\Big)\\
\le& \,C\cdot\|a\|^2_{L^2_w}.
\end{split}
\end{equation*}
We are done.
\end{proof}
We are now in a position to give the proof of Theorem 3.
\begin{proof}[Proof of Theorem 3]
As in the proof of Theorem 1, we write
\begin{equation*}
\begin{split}
&\sigma^p\cdot w(\{x\in\mathbb R^n:|g^*_{\lambda,\alpha}(a)(x)|>\sigma\})\\
\le&\sigma^p\cdot w(\{x\in Q^*:|g^*_{\lambda,\alpha}(a)(x)|>\sigma\})+\sigma^p\cdot w(\{x\in(Q^*)^c:|g^*_{\lambda,\alpha}(a)(x)|>\sigma\})\\
=&K_1+K_2.
\end{split}
\end{equation*}
Note that $\lambda>3+{(2\alpha)}/n>1$. Applying Chebyshev's inequality, H\"older's inequality, Lemma A and Lemma 4.1, we thus have
\begin{align}
K_1&\le\int_{Q^*}|g^*_{\lambda,\alpha}(a)(x)|^pw(x)\,dx\notag\\
&\le\left(\int_{Q^*}|g^*_{\lambda,\alpha}(a)(x)|^2w(x)\,dx\right)^{p/2}\left(\int_{Q^*}w(x)\,dx\right)^{1-p/2}\notag\\
&\le \|g^*_{\lambda,\alpha}(a)\|^p_{L^2_w}w(Q^*)^{1-p/2}\notag\\
&\le C\cdot\|a\|^p_{L^2_w}w(Q)^{1-p/2}\notag\\
&\le C\notag.
\end{align}
\par
We now turn to deal with $K_2$. In the proof of Theorem 2, we have already showed
\begin{equation}
|S_\alpha(a)(x)|^2\le C\left(\frac{1}{w(Q)^{1/p}}\right)^2.
\end{equation}
For any given $(y,t)\in\Gamma_{2^k}(x)$, $x\in(Q^*)^c$, then a simple calculation shows that $t\ge\frac{|x-x_0|}{2^{k+2}}$, $k\in\mathbb Z_+$. Hence, by the estimates (2) and (3), we can get
\begin{equation}
\begin{split}
|S_{\alpha,2^k}(a)(x)|^2&=\iint_{\Gamma_{2^k}(x)}\Big(\sup_{\varphi\in{\mathcal C}_\alpha}|a*\varphi_t(y)|\Big)^2\frac{dydt}{t^{n+1}}\\
&\le C\left(\frac{r^{n+\alpha}}{w(Q)^{1/p}}\right)^2\int_{\frac{|x-x_0|}{2^{k+2}}}^\infty\int_{|y-x|<2^kt}\frac{dydt}{t^{2n+2\alpha+n+1}}\\
&\le C\cdot2^{k(3n+2\alpha)}\left(\frac{r^{n+\alpha}}{w(Q)^{1/p}}\right)^2\frac{1}{|x-x_0|^{2n+2\alpha}}\\
&\le C\cdot2^{k(3n+2\alpha)}\left(\frac{1}{w(Q)^{1/p}}\right)^2.
\end{split}
\end{equation}
It follows immediately from the inequalities (5), (7) and (8) that
\begin{equation*}
\begin{split}
\big(g^*_{\lambda,\alpha}(a)(x)\big)^2&\le C\left(\frac{1}{w(Q)^{1/p}}\right)^2\bigg(1+\sum_{k=1}^\infty2^{-k\lambda n}\cdot2^{k(3n+2\alpha)}\bigg)\\
&\le C\left(\frac{1}{w(Q)^{1/p}}\right)^2,
\end{split}
\end{equation*}
where the last series is convergent since $\lambda>3+(2\alpha)/n$. Again, the rest of the proof is exactly the same as that of Theorem 1, we finally obtain
$$K_2\le C.$$
Therefore, we conclude the proof of Theorem 3.
\end{proof}

\end{document}